\documentclass[11pt,reqno]{amsart}
\usepackage{color}

\usepackage{amsfonts,amscd}
\usepackage{amssymb}
\usepackage{url}
\usepackage[english]{babel}

\usepackage{booktabs}

\theoremstyle{plain}
\newtheorem{theorem}                 {Theorem}      [section]

\theoremstyle{definition}

\newtheorem{definition}   [theorem]  {Definition}

\newtheorem{remark}       [theorem]  {Remark}

\numberwithin{equation}{section}

\def \cn{{\mathbb C}}
\def \hn{{\mathbb H}}

\def \rn{{\mathbb R}}

\def \E{\mathcal E}

\def\nab#1#2{\hbox{$\nabla$\kern -.3em\lower 1.0 ex
		\hbox{$#1$}\kern -.1 em {$#2$}}}

\def \Re{\mathfrak R\mathfrak e}

\def \u{\mathfrak u}

\def \GLR#1{\text{\bf GL}_{#1}(\rn)}

\def \GLC#1{\text{\bf GL}_{#1}(\cn)}
\def \glc#1{\mathfrak{gl}_{#1}(\cn)}
\def \GLH#1{\text{\bf GL}_{#1}(\hn)}

\def \SLR#1{\text{\bf SL}_{#1}(\rn)}
\def \tSL2{\widetilde{\text{\bf SL}}_{2}(\rn)}

\def \SLC#1{\text{\bf SL}_{#1}(\cn)}

\def \O#1{\text{\bf O}(#1)}
\def \SO#1{\text{\bf SO}(#1)}

\def \SOs#1{\text{\bf SO}^*(#1)}

\def \SOO#1#2{\text{\bf SO}(#1,#2)}

\def \U#1{\text{\bf U}(#1)}
\def \u#1{\mathfrak{u}(#1)}
\def \Us#1{\text{\bf U}^*(#1)}

\def \UU#1#2{\text{\bf U}(#1,#2)}

\def \SU#1{\text{\bf SU}(#1)}
\def \su#1{\mathfrak{su}(#1)}
\def \SUs#1{\text{\bf SU}^*(#1)}

\def \SUU#1#2{\text{\bf SU}(#1,#2)}

\def \Sp#1{\text{\bf Sp}(#1)}

\def \Spp#1#2{\text{\bf Sp}(#1,#2)}

\def \SpR#1{\text{\bf Sp}(#1,\rn)}

\DeclareMathOperator{\Div}{div} 

\DeclareMathOperator{\trace}{trace}

\numberwithin{equation}{section}
\allowdisplaybreaks

\begin{document}

\title[Proper $r$-harmonic functions from Riemannian manifolds]{Proper $r$-harmonic functions from\\ Riemannian manifolds}

\dedicatory{version 1.0040 - 24 November 2019}

\author{Sigmundur Gudmundsson}
\address{Mathematics, Faculty of Science\\
	University of Lund\\
	Box 118, Lund 221\\
	Sweden}
\email{Sigmundur.Gudmundsson@math.lu.se}

\author{Marko Sobak}
\address{Mathematics, Faculty of Science\\
	University of Lund\\
	Box 118, Lund 221\\
	Sweden}
\email{marko.sobak1@gmail.com}

\begin{abstract}
We introduce a new method for constructing complex-valued $r$-harmonic functions on Riemannian manifolds.  We then apply this for the important semisimple Lie groups $\SO n$, $\SU n$, $\Sp n$, $\SLR n$, $\SpR n$, $\SUU pq$, $\SOO pq$, $\Spp pq$, $\SOs{2n}$ and $\SUs{2n}$.
\end{abstract}



\maketitle

\section{Introduction}

Biharmonic functions are important in physics and arise in areas of continuum mechanics, including elasticity theory and the solution of Stokes flows. Here the literature is vast, but usually the domains are either surfaces or open subsets of flat Euclidean space, with only very few exceptions, see for example \cite{Bai-Far-Oua}. The development of the very last years has changed this and can be traced e.g. in the following publications: \cite{Gud-Mon-Rat-1}, \cite{Gud-13}, \cite{Gud-14}, \cite{Gud-15}, \cite{Gud-Sif-1}.  There the authors develop methods for constructing explicit $r$-harmonic functions on classical Lie groups and even some symmetric spaces.
\smallskip

In this paper we introduce a general scheme for constructing complex-valued $r$-harmonic functions $\phi:(M,g)\to\cn$ on Riemannian manifolds.  We then show that this can be employed in the important cases of the classical compact simple Lie groups $\SO n$, $\SU n$, $\Sp n$ and furthermore in the non-compact cases of  $\SLR n$, $\SpR n$, $\SUU pq$, $\SOO pq$, $\Spp pq$, $\SOs{2n}$ and $\SUs{2n}$.

\section{Preliminaries}\label{section-preliminaries}

Let $(M,g)$ be a smooth manifold equipped with a Riemannian metric $g$.  We complexify the tangent bundle $TM$ of $M$ to $T^{\cn}M$ and extend
the metric $g$ to a complex-bilinear form on $T^{\cn}M$.  Then the
gradient $\nabla \phi$ of a complex-valued function $\phi:(M,g)\to\cn$ is a
section of $T^{\cn}M$.  In this situation, the well-known linear
{\it Laplace-Beltrami} operator (alt. {\it tension} field) $\tau$ on $(M,g)$
acts locally on $\phi$ as follows
$$
\tau(\phi)=\Div (\nabla \phi)=\frac{1}{\sqrt{|g|}}
\frac{\partial}{\partial x_j}\left(g^{ij}\, \sqrt{|g|}\,
\frac{\partial \phi}{\partial x_i}\right).
$$
For two complex-valued functions $\phi,\psi:(M,g)\to\cn$ we have the
following well-known relation
\begin{equation*}
\tau(\phi\,\psi)=\tau(\phi)\,\psi+2\,\kappa(\phi,\psi)+\phi\,\tau(\psi),
\end{equation*}
where the {\it conformality} operator $\kappa$ is given by
$$
\kappa(\phi,\psi)=g(\nabla \phi,\nabla \psi).
$$

\smallskip

For a positive integer $r$, the iterated Laplace-Beltrami operator
$\tau^r$ is defined by
$$
\tau^{0} (\phi)=\phi,\quad \tau^r (\phi)=\tau(\tau^{(r-1)}(\phi)).
$$
\begin{definition}\label{definition-proper-r-harmonic} For a positive integer $r$, we say that a complex-valued function $\phi:(M,g)\to\cn$ is
\begin{enumerate}
\item[(a)] {\it $r$-harmonic} if $\tau^r (\phi)=0$,
\item[(b)] {\it proper $r$-harmonic} if $\tau^r (\phi)=0$ and
$\tau^{(r-1)}(\phi)$ does not vanish identically.
\end{enumerate}
\end{definition}

It should be noted that the {\it harmonic} functions are exactly the
$1$-harmonic and the {\it biharmonic} functions are the $2$-harmonic
ones. In some texts, the $r$-harmonic functions are also called {\it polyharmonic} of order $r$.
We also note that, if a function is $r$-harmonic, then it is also $p$-harmonic for any $p \geq r$.
Hence, one is usually interested in studying functions which are {\it proper} $r$-harmonic.

\begin{definition}\label{defi-eigenfamily}
Let $(M,g)$ be a Riemannian manifold.  Then a set
$$\E =\{\phi_i:M\to\cn\ |\ i\in I\}$$ of complex-valued functions is said to be an {\it eigenfamily} on $M$ if there exist complex
numbers $\lambda,\mu\in\cn$ such that 
$$\tau(\phi)=\lambda\,\phi\ \ \text{and}\ \ \kappa(\phi,\psi)=\mu\,\phi\,\psi$$ 
for all $\phi,\psi\in\E$.
\end{definition}

\begin{definition}\label{definition-eigenfunction}
Let $(M,g)$ be a Riemannian manifold.  Then a complex-valued function $\phi:M\to\cn$ is said to be an {\it eigenfunction} if it is eigen both with respect to the Laplace-Beltrami operator $\tau$ and the conformality operator $\kappa$ i.e. there exist complex numbers $\lambda,\mu\in\cn$ such that 
$$\tau(\phi)=\lambda\, \phi\ \ \text{and}\ \ \kappa(\phi,\phi)=\mu\, \phi^2.$$
\end{definition}

\begin{remark}\label{remark-eigenfunction-eigenfamily}
It is clear that any element $\phi\in\E$ of an eigenfamily is an eigenfunction in the sense of Definition \ref{definition-eigenfunction}.
\end{remark}

\section{A new construction method}
\label{section-construction}

Let $\phi:(M,g)\to\cn$ be a complex-valued function and $f:U\to\cn$ be a holomorphic function defined on an open subset $U$ of $\cn$ containing the image $\phi(M)$ of $\phi$. Then it is a direct consequence of the chain rule that the composition $\Phi=f\circ\phi:(M,g)\to\cn$ is harmonic if and only if  
\begin{equation}\label{general-equation}
\tau(\Phi)=\tau(f\circ\phi)=f''(\phi)\,\kappa(\phi,\phi) + f'(\phi)\,\tau(\phi)=0.
\end{equation}
In general this differential equation seems difficult, if not impossible, to solve.  But as we will now see, this problem can actually be solved in the case when $\phi:(M,g)\to\cn$ is an eigenfunction in the sense of Definition \ref{definition-eigenfunction}.  In that case the equation (\ref{general-equation}) takes the following form
\begin{equation*}
\tau(\Phi)=\tau(f\circ\phi)= \mu\,\phi^2\, f''(\phi) + \lambda\, \phi\, f'(\phi)=0.
\end{equation*}
Thus, $\Phi$ will be harmonic if and only if $f$ satisfies the complex ordinary differential equation
\begin{equation*}
\mu\,z^2\, f''(z) + \lambda\, z\, f'(z)=0.
\end{equation*}
Note that this equation has a singularity at $z = 0$, so we are interested in solving it on a simply connected subset of $\cn \setminus \{0\}$.
The canonical choice here is of course the slit plane $\cn \setminus (-\infty, 0]$.
The equation above can then be solved for $f : \cn \setminus (-\infty,0] \to \cn$ by elementary means.
Since $f$ is not defined on the entire complex plane, we are also forced to restrict the domain of
the resulting harmonic function $\Phi = f\circ \phi$ to the open set
\begin{equation*}
W=\{ x\in M \,\mid\, \phi(x) \not\in (-\infty,0] \} \subset M.
\end{equation*}
Having constructed a harmonic function, we can now use it to construct a biharmonic function.
Let $\Phi_2 = f_2 \circ \phi$ for some holomorphic function $f_2 : \cn \setminus (-\infty,0] \to \cn$.
Then we can make $\Phi_2$ proper biharmonic by requiring that it solves the Poisson equation
\begin{equation*}
\tau(\Phi_2) = \Phi. 
\end{equation*} 
Employing the chain rule, we see that this equation is equivalent to
\begin{equation*}
\mu\,\phi^2\, f_2''(\phi) + \lambda\, \phi\, f_2'(\phi) = f(\phi),
\end{equation*}
i.e.\ we want $f_2$ to satisfy the complex ordinary  differential equation
\begin{equation*}
\mu\,z^2\, f_2''(z) + \lambda\, z\, f_2'(z) = f(z), \quad z \in \cn \setminus (-\infty,0],
\end{equation*}
which can again be solved using elementary methods.
This procedure can be continued inductively to obtain  proper $r$-harmonic functions.
Indeed, suppose that a holomorphic function $f_{r-1}:\cn \setminus (-\infty,0]\to\cn$ has been constructed so that $\Phi_{r-1} = f_{r-1} \circ \phi$ is proper $(r-1)$-harmonic.
Then we can consider the function $\Phi_r = f_r \circ \phi : W\to\cn$, where $f_r: \cn \setminus (-\infty,0] \to \cn$ is holomorphic, and study the Poisson equation
\begin{equation*}
\tau(\Phi_{r}) = \Phi_{r-1}, 
\end{equation*}
which, as before, reduces to a complex ordinary differential equation that can easily be solved.

To state the solutions furnished by the method that we have just described, we will make extensive use of 
the holomorphic principal logarithm $\log:\cn \setminus (-\infty, 0]\to\cn$ and the standard notation
$$z^\alpha = \exp(\alpha \log(z)).
$$
As already hinted, our proper $r$-harmonic functions will be defined on the open subset $W$ of $M$ given by 
\begin{equation*}
W=\{ x\in M \,\mid\, \phi(x) \not\in (-\infty,0] \}.
\end{equation*}
We wish this domain to be as large as possible, so it might be convenient to introduce branch cuts other than $(-\infty,0]$, depending on the values of $\phi$.
For simplicity, we state our result on $\cn \setminus (-\infty,0]$.

\begin{theorem}
Let $\phi:(M,g)\to\cn$ be a complex-valued eigenfunction on a Riemannian manifold such that the tension field $\tau$ and the conformality operator $\kappa$ satisfy 
$$\tau(\phi)=\lambda\, \phi\ \ \text{and}\ \ \kappa(\phi,\phi)=\mu\,\phi^2$$
for some $\lambda, \mu \in \cn$.
Then for a natural number $r\ge 1$ and $(c_1,c_2)\in\cn^2$ any non-vanishing function 
$$\Phi_r:W=\{ x\in M \,\mid\, \phi(x) \not\in (-\infty,0] \}\to\cn$$ satisfying 
$$\Phi_r(x)= 
\begin{cases}
c_1\,\log(\phi(x))^{r-1}, 																& \text{if }\; \mu = 0, \; \lambda \not= 0\\[0.2cm]
c_1\,\log(\phi(x))^{2r-1}+ c_{2}\,\log(\phi(x))^{2r-2}, 								& \text{if }\; \mu \not= 0, \; \lambda = \mu\\[0.2cm]
c_{1}\,\phi(x)^{1-\frac\lambda{\mu}}\log(\phi(x))^{r-1} + c_{2}\,\log(\phi(x))^{r-1},	& \text{if }\; \mu \not= 0, \; \lambda \not= \mu
\end{cases}
$$ 
is proper $r$-harmonic on its open domain $W$ in $M$.
\end{theorem}

\medskip

\section{Eigenfunctions on semisimple Lie groups} \label{section-eigenfunctions-simple-semisimple}

In the last section we have introduced a new method for constucting $r$-harmonic functions on Riemannian manifolds in terms of eigenfunctions.  This method is not worth much without the existence of such objects.  The goal of this section is to show that these certainly do exist in many of the important cases of the classical semisimple Lie groups. 

As already pointed out in Remark \ref{remark-eigenfunction-eigenfamily} the elements $\phi\in\E$ of a given eigenfamily are automatically eigenfunctions i.e.\ exactly the objects that we are seeking. Eigenfamilies on classical semisimple Lie groups have already been constructed in \cite{Gud-Mon-Rat-1}, \cite{Gud-Sak-1} and \cite{Gud-Sak-2}.
We use the results from these papers to present interesting examples, see Table \ref{table-eigenfunctions}.

Let us discuss our notation and recall the definitions of some of the classical semisimple Lie groups.
We will denote the $n\times n$ identity matrix by $I_n$ and we will also use the standard notations
\begin{equation*}
J_n = \begin{bmatrix}
0 & I_n\\
-I_n & 0
\end{bmatrix}, \quad
I_{p,q} = \begin{bmatrix}
-I_p & 0\\
0 & I_q
\end{bmatrix}.
\end{equation*}

The Lie groups in this section are all assumed to be equipped with the left-invariant Riemannian metric induced by the canonical inner product
\begin{equation*}
g(Z,W) = \Re\trace(ZW^*)
\end{equation*}
on the Lie algebra $\glc{n}$.
We will consider the well-known special linear group
\begin{equation*}
\SLR{n} = \{ x \in \GLR{n} \,\mid\, \det x = 1 \}.
\end{equation*}
Moreover we will consider the important orthogonal, unitary, and quaternionic unitary Lie groups
\begin{eqnarray*}
	\O{n} &=& \{ x \in \GLR{n} \,\mid\, x \cdot x^t = I_n \},\\
	\U{n} &=& \{ z \in \GLC{n} \,\mid\, z \cdot z^* = I_n \},\\
	\Sp{n} &=& \{ q \in \GLH{n} \,\mid\, q \cdot q^* = I_n \},
\end{eqnarray*}
as well as the special orthogonal and the special unitary groups
\begin{eqnarray*}
	\SO{n} &=& \O{n} \cap \SLR{n},\\
	\SU{n} &=& \U{n} \cap \SLC{n}.
\end{eqnarray*}
We recall that their generalisations $\mathbf{O}(p,q),\UU{p}{q}, \Spp{p}{q}$ are defined as the groups of linear transformations on
$\rn^{p+q}$, $\cn^{p+q}$, $\hn^{p+q}$, respectively,  preserving a bilinear (or Hermitian) form induced by $I_{p,q}$.
More precisely,
\begin{eqnarray*}
	\mathbf{O}(p,q) &=& \{ x \in \GLR{p+q} \,\mid\, x \cdot I_{p,q} \cdot x^t = I_{p,q} \},\\
	\UU{p}{q} &=& \{ z \in \GLC{p+q} \,\mid\, z \cdot I_{p,q} \cdot z^* = I_{p,q} \},\\
	\Spp{p}{q} &=& \{ q \in \GLH{p+q} \,\mid\, q \cdot I_{p,q} \cdot q^* = I_{p,q} \},
\end{eqnarray*}
and we also put
\begin{eqnarray*}
	\SOO{p}{q} &=& \mathbf{O}(p,q) \cap \SLR{p+q},\\
	\SUU{p}{q} &=& \UU{p}{q} \cap \SLC{p+q}.
\end{eqnarray*}
For the quaternionic unitary group $\Spp{p}{q}$ we use the standard complex representation of $\GLH{p+q}$ in $\cn^{2(p+q)\times 2(p+q)}$ given by
$$
(z+jw)\mapsto q=\begin{bmatrix}z & w \\ -\bar w & \bar
z\end{bmatrix}.
$$
We will also consider the symplectic group
\begin{equation*}
\SpR{n} = \{ q \in \SLR{2n} \,|\, q \cdot J_n \cdot q^t = J_n \}.
\end{equation*}
Note that every element $q \in \SpR{n}$ has a unique representation as
\begin{equation*}
	q = \begin{bmatrix}
		x & y\\
		z & w
	\end{bmatrix},
\end{equation*}
for some matrices $x,y,z,w \in \rn^{n\times n}$.
Finally, we will consider the Lie groups
\begin{eqnarray*}
	\SOs{2n} &=& \{ q\in \SUU{n}{n} \,\mid\, q \cdot I_{n,n} \cdot J_n \cdot q^t = I_{n,n} \cdot J_n \},\\
	\SUs{2n} &=& \Us{2n} \cap \SLC{2n},
\end{eqnarray*}
where
\begin{equation*}
	\Us{2n} =
	\left\{
	z + jw =
	\begin{bmatrix}
	z & w\\
	-\overline{w} & \overline{z}
	\end{bmatrix}
	\,\mid\,
	z,w \in \GLC{n}
	\right\}.
\end{equation*}
Note that for any element $q \in \SOs{2n}$, we have a unique representation
\begin{equation*}
	q = \begin{bmatrix}
	z & w\\
	-\overline{w} & \overline{z}
\end{bmatrix},
\end{equation*}
where $z,w$ are $n\times n$ matrices.

Having recalled the definitions of these classical semisimple Lie groups, we are ready to present examples of eigenfunctions, see Table \ref{table-eigenfunctions}.
The reader should note that in Table \ref{table-eigenfunctions} we interpret the vectors as row vectors,
so that if $a,v \in \cn^n$, then $a^t v$ is an $n\times n$ matrix with entries $a_i v_j$.
We also denote
\begin{equation*}
\cn_1^p = \{ (z,w) \in \cn^p \times \cn^q \,\mid\, w = 0 \}, \quad
\cn_2^q = \{ (z,w) \in \cn^p \times \cn^q \,\mid\, z = 0 \}.
\end{equation*}
The machinery required for calculating the tension field and the conformality operator of the functions in Table \ref{table-eigenfunctions} 
can be found in \cite{Gud-Mon-Rat-1}, 
\cite{Gud-Sak-1} and \cite{Gud-Sak-2}, and is therefore omitted from this paper.
Most of the results presented here are immediate consequences of the calculations made there,
the only exception being the results concerning the Lie groups $\SLR{n}, \SU{n}, \SUU{p}{q}$ and $\SUs{2n}$,
which require a simple modification.

Let us describe this for the case of $\SU{n}$, the other cases follow the same pattern.
Following Table \ref{table-eigenfunctions}, define $\phi : \SU{n} \to \cn$ by
\begin{equation*}
\phi(z) = \trace(a^tvz^t),
\end{equation*}
where $a,v\in\cn^n$.
Note that this function extends naturally to $\U{n}$.
If we let $\hat{\tau}$ and $\hat{\kappa}$ denote the tension field and the conformality operator of $\U{n}$, respectively,
then as a consequence of Lemma 5.1 and Theorem 5.2 in \cite{Gud-Sak-1}, we get
\begin{equation*}
\hat{\tau}(\phi) = -n\,\phi \quad\text{and}\quad \hat{\kappa}(\phi,\phi) = -\phi^2.
\end{equation*}
Now it is easy to see that the matrix
\begin{equation*}
X = \frac{i}{\sqrt n} \, I_n \in \u{n}
\end{equation*}
is of unit length and generates the orthogonal complement of $\su{n}$ in $\u{n}$.
If we denote the tension field and the conformality operator of $\SU{n}$ by $\tau$ and $\kappa$ respectively, we therefore obtain
\begin{eqnarray*}
\tau(\phi) &=& \hat{\tau}(\phi) - X^2(\phi) = -n \, \phi + \frac{1}{n} \, \phi = -\frac{n^2-1}{n} \,\phi,\\[0.1cm]
\kappa(\phi,\phi) &=& \hat{\kappa}(\phi,\phi) - X(\phi)^2 = - \phi^2 + \frac{1}{n} \, \phi^2 = -\frac{n-1}{n}\,\phi,
\end{eqnarray*}
as claimed in Table \ref{table-eigenfunctions}.

\newpage

\vspace*{\fill}

\renewcommand{\arraystretch}{2}
\begin{table}[h]
\makebox[\textwidth][c]{
\begin{tabular}{ccccc}
\midrule\midrule
Lie group 	& Eigenfunction $\phi$  			& $\tau(\phi)$ 				& $\kappa(\phi)$ 			& Conditions
\\\midrule\midrule
$\SLR{n}$	& $\trace(Ax^t)$				& $\frac{n-1}{n}\,\phi$			& $-\frac{1}{n}\,\phi^2$	&   $A\in\cn^{n\times n}, \, AA^t = 0$
\\\midrule
$\SO{n}$	& $\trace(a^tvx^t)$			 	& $-\frac{n-1}{2}\,\phi$			& $-\frac{1}{2}\,\phi^2$ 	& $a,v \in\cn^n, \, v$ isotropic 
\\\midrule
$\SU{n}$	& $\trace(a^tvz^t)$				& $-\frac{n^2-1}{n}\,\phi$		& $-\frac{n-1}{n}\,\phi^2$& $a,v\in\cn^n$
\\\midrule
$\Sp{n}$	& $\trace(a^tvz^t + a^tuw^t)$	& $-\frac{2n+1}{2}\,\phi$		& $-\frac{1}{2}\,\phi^2$	& $a,v,u \in \cn^n$
\\\midrule
$\SOO{p}{q}$& $\trace(a^tvx^t)$			 	& $\frac{q-p+1}{2}\,\phi$		& $-\frac{1}{2}\,\phi^2$ 	& $a\in\cn^n, \, v \in \cn^p_1, \, v$ isotropic \\
			& $\trace(a^tvx^t)$			 	& $\frac{p-q+1}{2}\,\phi$		& $-\frac{1}{2}\,\phi^2$ 	& $a\in\cn^n, \, v \in \cn^q_2, \, v$ isotropic
\\\midrule
$\SUU{p}{q}$& $\trace(a^tvz^t)$				& $\frac{q^2-p^2+1}{p+q}\,\phi$	& $-\frac{p+q-1}{p+q}\,\phi^2$  & $a\in\cn^n, \, v \in \cn^p_1$ \\
			& $\trace(a^tvz^t)$			 	& $\frac{p^2-q^2+1}{p+q}\,\phi$	& $-\frac{p+q-1}{p+q}\,\phi^2$  & $a\in\cn^n, \, v \in \cn^q_2$
\\\midrule
$\Spp{p}{q}$& $\trace(a^tvz^t + a^tuw^t)$	& $-\frac{2(p-q)+1}{2}\,\phi$	& $-\frac{1}{2}\,\phi^2$		& $a \in \cn^n, \, v,u \in \cn^p_1$ \\
			& $\trace(a^tvz^t + a^tuw^t)$	& $-\frac{2(q-p)+1}{2}\,\phi$	& $-\frac{1}{2}\,\phi^2$		& $a \in \cn^n, \, v,u \in \cn^q_2$
\\\midrule
$\SpR{n}$	& $\trace( a^tv\,(x+iy)^t )$  	& $\frac{1}{2}\,\phi$ 			& $-\frac{1}{2}\,\phi^2$ 	& $a,v \in \cn^n$ \\
			& $\trace( a^tv\,(z+iw)^t )$  	& $\frac{1}{2}\,\phi$ 			& $-\frac{1}{2}\,\phi^2$ 	& $a,v \in \cn^n$
\\\midrule
$\SOs{2n}$	& $\trace(a^tvz^t)$				& $-\frac{1}{2}\,\phi$			& $-\frac{1}{2}\,\phi^2$  		& $a,v\in\cn^n$ \\
			& $\trace(a^tvw^t)$			 	& $-\frac{1}{2}\,\phi$			& $-\frac{1}{2}\,\phi^2$ 		& $a,v\in\cn^n$
\\\midrule
$\SUs{2n}$ 	& $\trace(Az^t + Bw^t)$ 		& $-\frac{2n+1}{2n}\,\phi$ 		& $-\frac{1}{2n}\,\phi^2$	& $A,B \in \cn^{n\times n}, \, AB^t = BA^t$
\\\midrule\midrule
\end{tabular}
}

\bigskip

\caption{Eigenfunctions on classical semisimple Lie groups.}
\label{table-eigenfunctions}

\end{table}

\vspace*{\fill}

\end{document}